\documentclass[12pt]{article}
\usepackage{eurosym}
\usepackage{amsmath, amsfonts, amssymb}

\def\CC{{\rm \kern.24em \vrule width.02em height1.4ex
depth-.05ex \kern-.26em C}}
\newtheorem{definition}{Definition}[section]
\newtheorem{theorem}[definition]{Theorem}
\newtheorem{lemma}[definition]{Lemma}
\newtheorem{proposition}[definition]{Proposition}
\newtheorem{corollary}[definition]{Corollary}

\newcommand{\be}{\begin{equation}}
\newcommand{\ee}{\end{equation}}
\newcommand{\bea}{\begin{eqnarray}}
\newcommand{\eea}{\end{eqnarray}}
\newcommand{\Bea}{\begin{eqnarray*}}
\newcommand{\Eea}{\end{eqnarray*}}

\begin{document}

\title{Ovoidal fibrations in $PG(3,q), q$ even}
\author{N.S. Narasimha Sastry and R.P. Shukla \\
Stat.- Math. Unit; Indian Statistical Institute; \\
8th Mile Mysore Road; R.V. College Post;\\
Bangalore - 560 059 (India) \\
and \\
Department of Mathematics; University of Allahabad; \\
Allahabad - 211 002 (India)}
\date{}
\maketitle

\abstract{We prove that, given a partition of the point-set of $PG(3,q), q=2^n >2$, by ovoids
$\{\theta_i\}^q_{i=0}$ of $PG(3,q)$ and a line $\ell$ of $PG(3,q)$, not tangent to
$\theta_0$ if $\ell^\perp$ denotes the polar of $\ell$ relative to the
symplectic form on $PG(3,q)$ whose isotropic lines are the tangent lines to
$\theta_0$, then $\ell$ and $\ell^\perp$ are tangent to distinct ovoids
$\theta_j, \theta_k$, both distinct from $\theta_0$ (Theorem \ref{main theorem}).
This uses the fact that the radical of the linear code generated by the dual duals
$\ell\cup \ell^\perp$ of the hyperbolic quadrics , with $\ell$ and $\ell^\perp$
 as above, is of codimension $1$ (Theorem \ref{thm 1.4}) }\newline

\section{Introduction and Statement of the Result}

An \emph{ovoid} in $PG(3,q)$ is a set of $q^{2}+1$ points, no three
collinear. Elliptic quadrics and the Suzuki-Tits ovoids (\cite{tits62}, See 
\cite{hirschfeld85}, 16.4 ), which exist, if and only if, $n$ is odd, are
the only known ovoids in $PG(3,q)$ and these are the only ovoids if $n\leq 6$
\cite{O'keefe96}. Classification of ovoids in $PG(3,q)$ is a fundamental
problem in Incidence Geometry. We mention in passing that the elliptic
quadrics are the only ovoid in $PG(3,q)$ if $q$ is odd (independently due to
Barlotti \cite{bar65} and Panella \cite{pan55}, see \cite{O'keefe96},
Theorem 2.1, p.178). We recall that a \textit{general linear complex} in $%
PG(3,q)$ is the set of all absolute lines with respect to a symplectic
polarity of $PG(3,q)$ (that is, a polarity for which all points are
absolute). If $P$ is the set of all points of $PG(3,q)$ and $L$ is a general
linear complex in $PG(3,q)$, then the incidence system $(P,L)$, which we
denote by $W(q)$, is a generalized quadrangle of order $(q,q)$ (see \cite%
{payne-thas84}, p.37). Since the collineation group $PGL(4,q)$ acts by
conjugation transitively on the set of all symplectic polarities of $PG(3,q)$%
, $W(q)$ is uniquely defined, up to a collineation of $PG(3,q)$. An ovoid of 
$W(q)$ is a set $\mathcal{O}$ of $q^{2}+1$ points, pairwise noncollinear in $%
W(q)$. Then, any line $l$ of $PG(3,q)$ meets $\mathcal{O}$ in a point if it
is a line of $W(q)$ and in zero or two points if it is not a line of $W(q)$
as shown by a simple count. Thus any ovoid of $W(q)$ is an ovoid of $PG(3,q)$%
. The following observation, due to Segre \cite{segre67}, proves the
converse: let $\mathcal{O}$ be an ovoid of $PG(3,q)$, $q$ even, and $W(%
\mathcal{O})$ be the incidence system whose points are the points of $%
PG(3,q) $ and lines are the tangent lines of $\mathcal{O}$. For each $x\in 
\mathcal{O}$, the union of the set of the $q+1$ tangent lines to $\mathcal{O}
$ through $x$ is a plane $\pi _{x}$. The correspondence $x\mapsto \pi _{x}$
defines a symplectic polarity (= null polarity) of $PG(3,q)$ \cite{dem64}
whose absolute lines are precisely the lines of $W(\mathcal{O})$. Thus, $W(%
\mathcal{O})=g(W(q))$ for some $g\in PGL(4,q)$ and the ovoid $\mathcal{O}$
of $PG(3,q)$ is an ovoid of (a unique copy of) $W(q)$. Hence, the
classification of ovoids of $PG(3,q)$, $q$ even, is equivalent to the
classification of ovoids in $W(q)$.

It is easy to see that any three mutually skew lines in $PG(3,q)$ have
exactly $q+1$ transversals (lines meeting each of the three given lines).
Such a set of $q+1$ transversals is called a \textit{regulus} in $PG(3,q)$.
The transversals to the lines in any regulus form another regulus, called
the \textit{opposite} of the given regulus. Thus, any three mutually skew
lines of $PG(3,q)$ are in a unique regulus. A \textit{spread} of $PG(3,q)$
is a set of $q^{2}+1$ mutually skew lines. A spread $\mathcal{S}$ is said to
be \textit{regular} (see \cite{dem64}, 5.1) if the unique regulus containing
three lines in $\mathcal{S}$ is contained in $\mathcal{S}$.

An \emph{ovoidal fibration} $\mathcal{F}=\{\theta _{i}\}_{i=0}^{q}$ of $%
PG(3,q)$ is a set of $q+1$ ovoids partitioning the point set of $PG(3,q)$.
We recall that a Singer group in $PGL(4,q)$ is a cyclic group of order $%
(q^{2}+1)(q+1)$ acting transitively on the point-set of $P(3,q)$. Let $S$ be
a Singer group in $PGL(4,q)$, and $T,K$ be its unique subgroups of orders $%
q^{2}+1$ and $q+1$, respectively. Then, the point $T$-orbits $%
\{E_{i}\}_{i=0}^{q}$ in $PG(3,q)$ is an ovoidal fibration by elliptic
quadrics (see \cite{bag-sas87}, Lemma 2, p.141); also \cite{Ebert85},
Theorem 3 and \cite{cossidente95}); the set $\mathcal{S}$ of common tangent
lines to the $E_{i}$'s is a regular spread in $PG(3,q)$; $T$ acts regularly
on the set $\mathcal{S}$ and $K$ acts regularly on each member of $\mathcal{S%
}$ \ (\cite{glynn88}; see \cite{bag-sas87}, Lemma 2, p.141). More generally,
if $\theta $ is an ovoid of $PG(3,q)$ and $\mathcal{S}$ is contained in the
general linear complex $\mathcal{L}=\mathcal{L}(\theta )$ consisting of all
tangent lines to $\theta $ (see, Segre \cite{segre67}) and if $K$ is the
group of collineations fixing each line in $\mathcal{S}$, then $K$ is a
subgroup of order $q+1$ contained in a Singer subgroup of $PGL(4,q)$ (see 
\cite{glynn88}) and the set $\{x(\theta ):x\in K\}$ of ovoids in $PG(3,q)$
is an ovoidal fibration of $PG(3,q)$ ( see \cite{Butler2008}, Theorem 3.1,
p.161). Further, $\theta $ is the only ovoid of the symplectic generalized
quadrangle $(P,\mathcal{L})\simeq W(q)$, where $P$ is the point-set of $%
PG(3,q)$ (see \cite{bag-sas87}, Lemma 2, p.141).

Let $\theta $ be an ovoid of $W(q)$ and `$\perp $' denote the orthogonality
relative to the nondengerate symplectic form on $PG(3,q)$ whose isotropic
lines are the tangent lines to $\theta $. Since the symplectic polarity
defined by $\theta $ maps an external line to $\theta $ to a secant line to $%
\theta $ and vice-versa, if $m$ is a line of $PG(3,q)$ which is not tangent
to $\theta $, then one of the lines $m,m^{\perp }$ shares two points with $%
\theta $ and the other is disjoint from it (\cite{hirschfeld85}, Corollary 1
of Theorem 16.1.8). Consequently, $\theta $ is an ovoid of $PG(3,q)$ also.
Thus, the classification of ovoids in $PG(3,q)$ and in $W(q)$ are equivalent
problems. We do not know any examples of ovoidal fibration containing a pair
of projectively nonequivalent ovoids (see (\cite{cossidente01}) for a
discussion of the case $q=8$); and of any examples of fibrations consisting
of projectively equivalent ovoids, but with no transitive action of a
subgroup of $PGL(4,q)$ on its constituent ovoids .

The subset $m\cup m^{\perp }$ of $P$, with $m,m^{\perp }$ as above, is
called a \emph{dual grid} of $W(q)$. With the lines of $W(q)$ incident with
it as `lines', $m\cup m^{\perp }$ is a $(1,q)$-subgeneralized quadrangle of $%
W(q)$. We denote by $\mathbb{D}$ the set of dual grids in $W(q)$.

\begin{proposition}
\label{Proposition 1}Let $\{\theta _{i}\}_{i=0}^{q}$ be an ovoidal fibration
of $PG(3,q)$ and $\mathcal{L}_{i}$ denote the general linear complex
consisting of all tangent lines to $\theta _{i},0\leq i\leq q$. Then, the
following hold:

\begin{itemize}
\item[(i)] The set $\mathcal{S}$ of common tangent lines to the ovoids $%
\theta _{i}$ is a regular spread in $PG(3,q)$. Further, $\{\mathcal{L}%
_{i}\}_{i=0}^{q}$ are the only general linear complexes in $PG(3,q)$
containing $\mathcal{S}$; and $\mathcal{L}_{i}\cap \mathcal{L}_{j}=\mathcal{S%
}$ for \ all $0\leq i\neq j\leq q$.

\item[(ii)] Each line $\ell $ of $PG(3,q)$ not in $\mathcal{S}$ is tangent
to a unique ovoid $\theta _{i}$, secant to $q/2$ ovoids $\theta _{j}$, $%
j\neq i$, and is disjoint from each of the remaining $q/2$ ovoids $\theta
_{k}$ of the fibration.
\end{itemize}
\end{proposition}

Our object in this note is to prove

\begin{theorem}
\label{main theorem}Let $\{\theta _{i}\}_{i=0}^{q}$ be an ovoidal fibration
of $PG(3,q)$, $W(q)$ be the generalized quadrangle whose line set is the set 
$\mathcal{L}(\theta _{0})$ of tangent lines to $\theta _{0}$ and $m$ be a
line of $PG(3,q)$ not in $\mathcal{L}(\theta _{0})$. Then, $m$ and its perp $%
m^{\perp }$ in $W(q)$ are tangent to distinct ovoids $\theta _{i}$, each
distinct from $\theta _{0}$.
\end{theorem}

An intrinsic description of the ovoids $\theta _{i}$ the lines $m$ and $%
m^{\perp }$ are tangent to may be interesting. We would like to view this
note as a contribution towards understanding the packings of $PG(3,q)$ by
ovoids.

Our proof of Theorem \ref{main theorem} uses a property of the binary code
we now define. Let $G$ denote the projective symplectic subgroup $PSp(4,q)$
\ of $PGL(4,q)$ defined by $W(q)$. Let $\mathbb{F}_{2}^{P}$ and $\mathbb{F}%
_{2}^{\mathbb{D}}$ denote the $\mathbb{F}_{2}G$-permutation modules on $P$
and $\mathbb{D}$, respectively, and $\eta $ denote the $\mathbb{F}_{2}G$%
-module homomorphism taking $m\cup m^{\perp }\in \mathbb{D}$ to $\sum_{x\in
m\cup m^{\perp }}x\in \mathbb{F}_{2}^{P}$. We identify the characteristic
function of a subset of $P$, considered as an element of $\mathbb{F}_{2}^{P}$%
, with the subset itself. Let $\mathcal{C}$ and $\mathcal{D}$ denote the $%
\mathbb{F}_{2}G$-submodules of $\mathbb{F}_{2}^{P}$ whose generators are,
respectively, the lines of $W(q)$ and the dual grids of $W(q)$. Then, $\eta (%
\mathbb{F}_{2}^{\mathbb{D}})=\mathcal{D}$. Since any line of $W(q)$ meets a
dual grid in zero or two points, $\mathcal{D}$ is contained in the dual code 
$\mathcal{C}^{\perp }$ of $\mathcal{C}$. Further, it is generated by words
of $\mathcal{C}^{\perp }$ of minimum weight (\cite{bag-sas88}, Theorem 1.4
). \smallskip

\begin{corollary}
\label{Corollary 3}$\mathcal{D}\underset{+}{\subset }\mathcal{C}^{\perp }$.
\end{corollary}

We need the following

\begin{theorem}
\label{thm 1.4}The $\mathbb{F}_{2}G$-radical $U$ of $\mathcal{D}$ is of
codimension one in $\mathcal{D}$. Consequently, the sum of two dual grids of 
$W(q)$ is in the radical $U$.
\end{theorem}

\begin{lemma}
\label{Lemma5}The group $G$ contains a unique conjugacy class of subgroups $%
T $ of order $q^{2}+1$ and $T$ is cyclic. Let $\{E_{i}\}_{i=0}^{q}$ be the
point $T$-orbits in $PG(3,q)$ and $\mathcal{S}$ denote the set of all common
tangent lines to $E_{i}$'s. Let $\ell $ be a line of $PG(3,q)$. Then, $%
\sum\limits_{t\in T}t(\ell )\in \mathbb{F}_{2}^{P}$ is $P$ (i.e., the
`all-one' vector) or the unique ovoid $E_{i}$ the line $\ell $ is tangent
to, according as $\ell \in \mathcal{S}$ or $\ell \not\in \mathcal{S}$.
\end{lemma}

\section{Preliminaries}

Let $V(q)$ be the vector space of dimension four over $\mathbb{F}_{q}$; $%
\hbar $ be a non-degenerate symplectic bilinear form on it; and $W(q)$ be
the incidence system with the set $P=PG(3,q)$ of all one dimensional
subspaces of $V(q)$ as its \textit{point-set}, the set $L$ of all two
dimensional subspaces of $V(q)$ which are isotropic with respect to $\hbar $
as its \textit{line-set} and symmetrized inclusion as the \textit{incidence}%
. Then, $W(q)$ is a regular generalized quadrangle of order $q$ (\cite%
{payne-thas84}, p.37) and the symplectic group $G$ defined by $\hbar $ acts
as incidence preserving permutations on the sets $P$ and $L$.\smallskip

Let $k$ be an algebraically closed extension field of $\mathbb{F}_{q}$. Let $%
N=\{0,1,\cdots ,2n-1\}$. Addition in $N$ is always taken modulo $2n$. Let $%
\mathcal{N}$ denote the set of all subsets $I$ of $N$ containing no
consecutive elements. Let $V=V\left( q\right) \otimes k$. The natural
extension of the symplectic form $\hbar $ to $V$ defined\ above is also
denoted by $\hbar $. Then, $G$ is the subgroup of the algebraic group $%
Sp(V)\simeq Sp(4,k)$ fixed by the $n$-th power of the Frobenius map $\sigma $
(which is the algebraic group endomorphism of $GL(4,k)$ raising each entry
of a matrix to its $2^{nd}$-power). It is well known that $Sp(V)$ has an
algebraic group endomorphism $\tau $ with $\tau ^{2}=\sigma $ (\cite%
{steinberg68}, Theorem 28, p.146). For any non-negative integer $i$, we
denote by $V_{i}$ the $Sp(V)$- module whose $k$- vector space structure is
the same as that of $V$ and an element $g$ of $Sp(V)$ acts on $V_{i}$ as $%
\tau ^{i}(g)$ would act on $V$. For $I\subseteq N $, \ let $V_{I}$ denote
the $kG$- module ${\otimes _{i\in I}V}_{i}$ (with $V_{\emptyset }=k$). Then,
by Steinberg's tensor product Theorem (\cite{steinberg63}, \S 11), $%
\{V_{I}:I\subseteq N\}$ is a complete set of inequivalent simple $kG$
-modules. For a $kG$-module $M$, we denote by $[M:V_{I}]$ the multiplicity
of $V_{I}$ in a composition series of $M$. We denote by rad$(M)$ the \textit{%
radical of $M$} (that is, the smallest submodule of $M$ with semisimple
quotient). We refer to $M/rad(M)$ as the \textit{head} of $M$.

We now describe a graph automorphism $\tau $ of $Sp(V),$ following (\cite%
{flesner75}, pp.58-60). (The argument presented in loc. cit. constructs a
graph automorphism for $G=Sp(V(q))$, however the arguments are valid for $%
Sp(V)$ also.) Let $\{e_{1}, e_{2}, e_{3},e_{4}\}$ be an ordered basis of $V$ 
and $Q$ denote the
nondegenerate quadratic form on the exterior square $\Lambda ^{2}V$ of $V$
defined by 
\begin{equation*}
Q\left( \Sigma _{1\leq i<j\leq 4}\lambda _{ij}e_{i}\wedge e_{j}\right)
=\lambda _{12}\lambda _{34}+\lambda _{13}\lambda _{24}+\lambda _{14}\lambda
_{23}
\end{equation*}%
(whose zero set in $P(\Lambda ^{2}V)$ is the well-known \emph{Klein quadric}%
). Let $\beta $ denote the polarization of $Q$ and $\gamma =e_{1}\wedge
e_{4}+e_{2}\wedge e_{3}\in \Lambda ^{2}V$. Then, the restriction of $Q$ to
the hyperplane $U=\{x\in \Lambda ^{2}V:\beta \left( x,\gamma \right) =0\}$
of $\Lambda ^{2}V$ is a nondegenerate quadratic form; and the restriction of 
$\beta $ to $U$ is an alternating form with radical $k\gamma .$ The
alternating form $\overline{\beta }$ induced by $\beta $ on $\overline{U}=%
\frac{U}{k\gamma }$ is nondegenerate. So the symplectic space $\left( 
\overline{U},\overline{\beta }\right) $ is isometric to $\left( V,\hbar
\right) .$ Let $\overline{p}:\overline{U}\rightarrow V$ be the isometric
isomorphism induced by the linear map $p:U\rightarrow V$ defined by 
\begin{align*}
p\left( e_{1}\wedge e_{2}\right) & =e_{1},\text{ }p\left( e_{1}\wedge
e_{3}\right) =e_{2},\text{ }p\left( e_{2}\wedge e_{4}\right) =e_{3}, \\
p\left( e_{3}\wedge e_{4}\right) & =e_{4},\text{ }p\left( \gamma \right) =0.
\end{align*}

\noindent Then the map taking $g\in Sp(V)$ to $\overline{p}(\overline{\wedge
^{2}(g))\text{ }}\overline{p}^{-1}\in Sp(V)$ is a graph automorphism $\tau $
of $Sp(V)$ which, on restriction to $G$, gives a graph automorphism of $G$.
\bigskip

We use the following results: \smallskip

\begin{lemma}
\label{c-r} (\cite{cu-reiner81}, Corollary 7.11, p.148) Let $K$ be a field
of characteristic $p>0$, let $E$ be a field extension of $K$ and let $X$ be
a finite group. Then, for each simple $KX$ -module $M,M^{E}=M\otimes _{K}E$
is a direct sum of simple $EX$ -modules, no two of which are isomorphic.
\end{lemma}

\begin{lemma}
\label{s-s2007} (\cite{s-s2007}, Corollary 6) Let $K,K^{\prime }\in \mathcal{%
N}$ be distinct and $f_{\hbar }$ be a nondegenerate quadratic form on $V(q)$
of index $2$ which polarizes to $\hbar $. Then, $V_{K}$ and $V_{K^{\prime }}$
are semi-simple $k\Omega (f_{\hbar })$-module with no irreducible factors in
common.
\end{lemma}

\section{Proofs}

\textbf{Proof of Theorem \ref{thm 1.4}:} Let $g\in \mathbb{D}$ and $L$ = Stab%
$_{G}(g)$. We view ${\mathbb{F}}_{2}^{\mathbb{D}}$ as the induced module of
the trivial $\mathbb{F}_{2}G$- module $\mathbb{F}_{2}$, and write Ind$%
_{L}^{G}(\mathbb{F}_{2})={\mathbb{F}}_{2}^{\mathbb{D}}$. Then, 
\begin{equation*}
\lbrack {\mathbb{F}}_{2}^{\mathbb{D}}/Rad({\mathbb{F}}_{2}^{\mathbb{D}}):%
\mathbb{F}_{2}]=1
\end{equation*}%
because, by Frobenius reciprocity (\cite{lang}, p.689), 
\begin{equation*}
dim_{\mathbb{F}_{2}}(Hom_{\mathbb{F}_{2}G}(Ind_{L}^{G}(\mathbb{F}_{2})),%
\mathbb{F}_{2})=dim_{\mathbb{F}_{2}}(Hom_{\mathbb{F}_{2}L}(\mathbb{F}_{2},%
\mathbb{F}_{2}))=1.
\end{equation*}%
This shows that the trivial module $\mathbb{F}_{2}$ is a summand of the head
(that is, the largest semi-simple quotient) of the $\mathbb{F}_{2}G$-module $%
{\mathbb{F}}_{2}^{\mathbb{D}}$. Let $\mathcal{H}$ denote the set of all
hyperbolic quadrics of $W(q)$ and $f_{\hbar }$ be a nondegenerate quadratic
form on $V(q)$ of index $2$ which polarizes to $\hbar $. Then, the variety
defined by $f_{\hbar }$ is a member of $\mathcal{H}$. Let $O(f_{\hbar
})\subseteq G$ denote the orthogonal group of $f_{\hbar }$ and $\Omega
(f_{\hbar })$ denote its commutator subgroup. Let $I$ be a nonempty subset
of $N$. Then by Lemma \ref{s-s2007}, 
\begin{equation*}
Hom_{k\Omega (f_{\hbar })}(k,V_{I})\cong Hom_{k\Omega (f_{\hbar
})}(V_{I_{e}},V_{I_{o}})=\{0\}.
\end{equation*}%
This proves that $V_{I}$ has no nonzero fixed points for the action of $%
O(f_{\hbar })$.

Since $G$ acts transitively on $\mathcal{H}$ and since the automorphism $%
\tau $ maps the stabilizer of a hyperbolic quadric in $W(q)$ to the
stabilizer of a dual grid in $W(q)$ and vice-versa, 
\begin{equation*}
dim_{k}(Hom_{kG}(Ind_{L}^{G}(k),V_{I}))=dim_{k}(Hom_{kG}(Ind_{O(f_{\hbar
})}^{G}(k),V_{I})).
\end{equation*}%
Again using Frobenius reciprocity, we have 
\begin{align*}
dim_{k}(Hom_{kG}(Ind_{L}^{G}(k),V_{I}))& =dim_{k}(Hom_{kG}(Ind_{O(f_{\hbar
})}^{G}(k,V_{I})) \\
& =dim_{k}(Hom_{kO(f_{\hbar })}(k,V_{I})) \\
& \leq dim_{k}(Hom_{k\Omega (f_{\hbar })}(k,V_{I}))=0
\end{align*}

This proves that the head of the $kG$-module $k^{\mathbb{D}}\cong
Ind_{L}^{G}(k)$ is $k$. Now the Lemma \ref{c-r} completes the proof of
Theorem \ref{thm 1.4}.

\textbf{Proof of Proposition \ref{Proposition 1}}: The statement (i) follows
from (\cite{bag-sas13}, Theorem 2.2) and (\cite{Butler2008}, Lemma 2.4 ).
Since $|\mathcal{L}_{i}|=(q^{2}+1)(q+1)$ and $\mathcal{L}_{i}\cap \mathcal{L}%
_{j}=\mathcal{S}$ which has $q^{2}+1$ elements for all $i,j,$ $i\neq j$ (see 
\cite{bag-sas13}, Corollary 3.3), each line $\ell $ of $PG(3,q)$ not in $%
\mathcal{S}$ is in $\mathcal{L}_{i}$ for a unique $i$ and $|\ell \cap \theta
_{j}|\in \{0,2\}$ for all $j\neq i$. Since $\{\theta _{i}\}$ partitions $%
PG(3,q)$, (ii) follows.

\textbf{Proof of Lemma \ref{Lemma5}:} The first statement follows from (\cite%
{flesner75}, Lemma 3); the second from (\cite{bag-sas87}, Lemma 2, p.141);
and the third follows by Proposition \ref{Proposition 1}(ii) and the regular
action of $T$ on each $E_{i}$ (see \cite{bag-sas87}, Lemma 4, p.142).

\textbf{Proof of Theorem \ref{main theorem}:} Let $\mathcal{S}$ be the set
of common tangent lines to the $\theta _{i}$'s. Then, $\mathcal{S}$ is a
regular spread in $PG(3,q)$ by Proposition \ref{Proposition 1} (i). Let $H$
be the subgroup of $PGL(4,q)$ fixing each line in $\mathcal{S}$. Then $H$ is
contained in a Singer subgroup $K$ of $PGL(4,q)$ (see \cite{glynn88}).
Consider the unique subgroup $T$ of $K$ of order $q^{2}+1$. The point-orbits
of $T$ are elliptic ovoids $E_{i}$, forming an ovoidal fibration of $PG(3,q)$
and $\mathcal{S}$ is also the set of common tangents to $E_{i}$'s (See \cite%
{bag-sas87}, Lemma 2, p.141). Let $\mathcal{L}_{i}$ be the set of all
tangent lines to $\theta _{i}$ and $W(q)=W(\mathcal{L}_{0})$. Since $%
\mathcal{S}$ is contained in exactly $q+1$ general linear complexes (see
Proposition \ref{Proposition 1} (i)), we may assume that $\mathcal{L}_{i}$
is also the set of all tangent lines to $E_{i}$. Thus $\theta _{0}$ and $%
E_{0}$ are ovoids of $W(q)=W(L_{0})$ and $T$ is a subgroup of $G$ (= $%
Sp(4,q) $ defined by $W(q)$). Now, assume that $m\cup m^{\perp }$ is a dual
grid of $W(q)$ such that $m$ amd $m^{\perp }$ are tangents to $\theta _{i}$
for some $i\neq 0$. Then, both are tangents to $E_{i}$ also. Since $%
m,m^{\perp }\notin \mathcal{S}$, by Lemma \ref{Lemma5}, 
\begin{equation*}
\sum_{t\in T}t(m+m^{\perp })=E_{i}+E_{i}=0.
\end{equation*}%
On the other hand, since for each $t\in T$, $t(m+m^{\perp })\in \mathcal{D}$
and the codimension of $U$ in $\mathcal{D}$ is one, $\sum_{t\in
T}t(m+m^{\perp })=m+m^{\perp }$ (mod $U$). Since $G$ is transitive on $%
\mathbb{D}$, $U$ can not contain any dual grid of W(q). This completes the
proof of the theorem. \smallskip

\textbf{Proof of Corollary \ref{Corollary 3}: }Let $T$ be a cyclic subgroup
of order $q^{2}+1$ contained in the stabilizer in $G$ of an elliptic ovoid $%
\theta _{0}$ of $W(q)$. Let $\{\theta _{i}\}_{i=0}^{q}$ be the point $T$%
-orbits in $PG(3,q)$ and $\mathcal{S}$ denote the set of all common tangent
lines to $\theta _{i}$'s (See \cite{bag-sas87}, Lemma 2, p.141). Consider
the $\mathbb{F}_{2}$-linear map $\sigma :\mathbb{F}_{2}^{P}\longrightarrow 
\mathbb{F}_{2}^{P}$ defined by $\sigma (w)=\sum\limits_{t\in T}t(w)$. As we
noted earlier, $\mathcal{D}\subseteq \mathcal{C}^{\perp }$. Assume that $%
\mathcal{C}^{\perp }=\mathcal{D}$. Then $\mathcal{C}=(\mathcal{C}^{\perp
})^{\perp }=\mathcal{D}^{\perp }\supseteq \mathcal{D}$, as any two dual
grids of $W(q)$ intersect in zero or two points. Now, $\sigma (\mathcal{C}%
)=\{\emptyset ,P,\theta _{0},P\setminus \theta _{0}\}$ (see Lemma \ref%
{Lemma5}); however if $m\cup m^{\perp }$ is a dual grid of $W(q)$, then by
Lemma \ref{Lemma5} and Theorem \ref{main theorem} $\sigma (m+m^{\perp
})=\theta _{i}+\theta _{j}$ for distinct $i$ and $j$ ($i>0,j>0$), a
contradiction.


\begin{thebibliography}{99}
\bibitem{bag-sas87} B. Bagchi and N.S.N. Sastry, Even order inversive
planes, generalized quadrangles and Codes, Geom. Dedicata, 22 (1987) 137-147

\bibitem{bag-sas88} B. Bagchi and N.S.N. Sastry, One-step completely
orthogonalizable codes from generalized quadrangles, Inform. and Comput. 77
(1988) 123-130

\bibitem{bag-sas13} B. Bagchi and N.S.N. Sastry, Ovoidal packing of $P(3,q)$
for even $q$, Discrete Math. 313 (2013) 2213-2217

\bibitem{bar65} A. Barlotti, Un'estensione del teorema di
Segre-Kustaanheimo, Boll. Un. Mat. ltal. 10 (1955) 96-98.

\bibitem{Butler2008} D.K.Butler, On the intersection of ovoids sharing a
polarity, Geom. Dedicata 135 (2008) 157-165

\bibitem{cossidente01} A. Cossidente and Sam K.J. Vereecke, Some geometry of
the isomorphism $Sp(4,q)\cong O(5,q)$, $q$ even, J. Geom. 70 (2001) 28-37

\bibitem{cossidente95} A. Cossidente and L. Storme, Caps on elliptic
quadrics, Finite Fields Appl. (1995) 412-420

\bibitem{cu-reiner81} C.W. Curtis and I. Reiner, \textit{Methods of
Representation Theory, Vol I}, Wiley-Interscience Pub. New York (1981)

\bibitem{dem64} P. Demboswki, \textit{Finite Geometries}, Springer-Verlag,
Berlin, Heidelberg, New York, (1967)


\bibitem{Ebert85} G.L. Ebert, Partitioning projective geometries into caps,
Canad. J.Math.37 (1985) 163-175

\bibitem{flesner75} D.E. Flesner, Maximal subgroups of PSp$_4(2^n)$
containing central elations or noncentered skew elatins, Illinois J. Math.
19 (1975), 247-268

\bibitem{glynn88} D.G. Glynn, On the set of lines of $PG(3,q)$ corresponding
to a maximal cap contained in the Klein quadric of $PG(5,q)$, Geom. Dedicata
26 (1988) 273-280

\bibitem{hirschfeld85} J.W.P. Hirschfeld, \textit{Projective spaces of three
dimensions}, Clarendon Press, Oxford (1985)

\bibitem{hirschfeld91} J.W.P. Hirschfeld, Projective spaces of square size,
Simon Stevin 65 (1991) 319-329

\bibitem{huppert-1} B. Huppert, \textit{Enddliche Gruppen },
Springer-Verlag, Berlin, Heidelberg, New York, (1967)

\bibitem{lang} S. Lang, \textit{Algebra}, Third Edition, Addison-Wessley,
1999.

\bibitem{O'keefe96} C.M. O'Keefe, Ovoids in $PG(3,q)$, a survey, Discrete
Math. 151(1996) 171-188

\bibitem{pan55} G.Panella, Caratterizzazione delle quadriche di uno spazio
(tridimensionale) lineare sopra un corpo finito. Boll. Un. Mat. Ital. 10
(1955) 507-513

\bibitem{payne-thas84} S.E. Payne and J.A. Thas, \textit{Finite generalized
quadrangles,} Advance Publishing Program, Pitman, Boston 1984

\bibitem{s-s2007} N.S.N. Sastry and R. P. Shukla, Structure of a code
related to Sp(4,q), q even, Proc. Indian Acad. Sci.(Math. Sci.) 117(4)
(2007), 457-470

\bibitem{segre67} B. Segre, \textit{Introduction to Galois Geometries,} Mem.
Acad. Naz. Lincei 8 (1967), 137-236

\bibitem{steinberg63} R. Steinberg, Representations of Algebraic groups,
Nagoya J. Math. 22 (1963) 33-56

\bibitem{steinberg68} R. Steinberg, \textit{Lectures on Chevalley Groups},
Mimeographed Notes, Yale Univ. Math. dept., New Haven, Conn., 1968

\bibitem{tits62} J. Tits, Ovoides et groupes de Suzuki, Arch. math.
13(1962), 187-198
\end{thebibliography}
\end{document}